\documentclass[12pt,oneside]{amsart}
\usepackage{geometry}                
\usepackage{graphicx}
\usepackage{amssymb}
\usepackage{epstopdf}

\usepackage{amsmath,amsthm,amscd,amssymb, enumerate}
\usepackage{latexsym}
\usepackage[colorlinks,citecolor=red,pagebackref,hypertexnames=false]{hyperref}
\numberwithin{equation}{section}
\theoremstyle{plain}
\newtheorem{theorem}{Theorem}[section]

\newtheorem{corollary}[theorem]{Corollary}

\theoremstyle{definition}

\theoremstyle{remark}

\newtheorem{case[theorem]}{Case}

\def\bc{\begin{corollary}}
\def\ec{\end{corollary}}
\def\be{\begin{equation}}
\def\ee{\end{equation}}

\def\ss{\smallskip}
\def\ms{\medskip}

\def\ni{\noindent}

\def\R{\mathbb R}
\def\hd{{\dim_{\mathcal H}}}

\def\R{\Bbb R}

\def\vt{\vec{t}\,}

\def\hd{\dim_{\mathcal H}}

\def\H{\mathbb H}

\def\({\left(}
\def\){\right)}
\def\[{\left[}
\def\]{\right]}
\def\<{\left\langle}
\def\>{\right\rangle}

\def\d{\partial}

\def\rt{\mathcal R_{t}}
\def\zvt{Z_{\vt}}
\def\rvt{\mathcal R_{\vt}}
\def\sone{\mathbb S^1}

\def\sd{\mathbb S^{d-1}}

\def\tsd{T\, \sd}
\def\sd{\mathbb S^{d-1}}
\def\om{\omega}

\def\g{\frak g}

\def\a{\alpha}
\def\b{\beta}
\def\ap{\a_{\Phi}}
\def\bp{\b_{\Phi}}
\def\az{\a_{\Phi}}

\def\z{\mathbf 0}
\def\mzs{\setminus\z}

\def\corank{\hbox{corank}}
\def\dim{\hbox{dim}}
\def\dist{\hbox{dist}}
\def\rank{\hbox{rank}}

\title[Configuration sets with nonempty interior]{Configuration sets with nonempty interior}

\date{}    

\author{Allan Greenleaf, Alex Iosevich and Krystal Taylor}

\email{allan@math.rochester.edu}

\email{iosevich@math.rochester.edu}

\email{taylor.2952@osu.edu}

\address{Department of Mathematics, University of Rochester, Rochester, NY 14627}

\address{Department of Mathematics, University of Rochester, Rochester, NY 14627}

\address{Department of Mathematics, Ohio State University, Columbus, Ohio 43210}

\begin{document}
\maketitle

\begin{abstract} A  theorem of Steinhaus states that if $E\subset \R^d$ has positive Lebesgue measure, then
the difference set $E-E$ contains a neighborhood of $0$.
Similarly, if $E$ merely has Hausdorff dimension $\hd(E)>(d+1)/2$,
a result of Mattila and Sj\"olin states that 
 the distance set $\Delta(E)\subset\R$ contains an open interval.
In this work, we study such results from a general viewpoint, 
replacing $E-E$ or $\Delta(E)$
 with more general  $\Phi\,$-{configurations} 
for a  class of $\Phi:\R^d\times\R^d\to\R^k$,
and showing that, under suitable lower bounds
on  $\hd(E)$  
and a regularity assumption on the family of generalized Radon transforms associated
 with $\Phi$, it follows that
the set $\Delta_\Phi(E)$ of  $\Phi$-configurations in $E$  has nonempty interior in $\R^k$. 
Further extensions hold for  $\Phi\,$-configurations generated by two sets, $E$ and $F$, in spaces of possibly different 
dimensions and   with suitable lower bounds on $\hd(E)+\hd(F)$.

\end{abstract}  

\maketitle


\section{Introduction} \label{sec intro}

A classical theorem of Steinhaus \cite{Steinhaus} states that if $E\subset \R^d, d\ge 1,$ with positive Lebesgue measure, $|E|_d>0$,
then the difference set $E-E\subset \R^d$ contains a neighborhood of the origin.
 $E-E$ can interpreted as the set of two-point configurations, $x-y$, of points of $E$ modulo the translation group.
 A variant of this  was obtained by Mattila and Sj\"olin \cite{MS99}
 for \emph{thin sets}, i.e.,   $E$ with $|E|_d=0$ 
but satisfying a lower bound on the Hausdorff dimension, $\hd(E)$,
in the context of the Falconer distance problem: 
 if $\Delta(E)$ is  the distance set of $E$,
 $\Delta(E):=\{|x-y|: x,y\in E\}\subset\R$,  
then if $\hd(E)>\frac{d+1}2$, it follows that $\Delta(E)$ contains an open interval.
The purpose of the current paper is to generalize these results in two ways: 
to  two-point configurations in $E$ as measured by a general class of $\Phi$\emph{-configurations}, 
which can be nontranslation-invariant, and indeed not even  in Euclidean space,
and  to allow  \emph{asymmetric} configurations between sets in different spaces, e.g.,  
between points and lines or points and circles in $\R^2$, or lines and lines in $\R^3$. 
In the process, we shall establish non-empty interior results for some configuration sets
for which previously it was not even known  that the configuration space has positive Lebesgue measure. 
\medskip

In order to formulate these more general results, consider the models  $E-E$ and $\Delta(E)$ as the images of $E\times E$ 
under the maps $(x,y)\to x-y$ and $(x,y)\to |x-y|$, resp. 
Now consider a $C^\infty$ function  $\Phi:\R^d\times\R^d\to R^k,\,  k\le d,$  
which is a  defining function (vector-valued if $k>1$) in the sense that the differential $D\Phi(x,y)$ has maximal rank everywhere. 
Thus, $\Phi$ is a submersion and hence for each $\vt\in\R^k$, the level set
\be\label{def zt}
Z_{\vt}:=\big\{(x,y)\in \R^d\times \R^d: \Phi(x,y)=\vt\big\}
\ee
is a smooth, codimension $k$ surface in $ \R^d\times \R^d$, 
and  the $\zvt$  form a family of incidence relations on $\R^d$, indexed by $\vt$. 
(For $k=1$,  the scalar will be denoted by $t$.)
More generally, since many $\Phi$ of interest, such as those defining the generalized distance sets discussed below in
Sec. \ref{sec examples},
have points in the domain where they fail to be smooth, or have critical values in the codomain,
it is useful to restrict the domain or codomain and consider $\Phi:X\times Y \to W$, where (for now) $X,\, Y\subset \R^d$ and $W\subset\subset R^k$. 
\smallskip

If  the assumption that $\rank\!\left(D\Phi\left(x,y\right)\right)=k$  is strengthened slightly to  the condition 
(which is standard in integral geometry) that 
\be\label{cond dPhi}
\rank(D_x\Phi)=k\,\hbox{ and }\rank(D_y\Phi)=k
\ee 
everywhere,  then each of the two projections,
$$\pi_X,\, \pi_Y:\zvt\to\R^d,\quad \pi_X(x,y):=x,\quad \pi_Y(x,y):=y,$$
are submersions, and $\zvt$ is a \emph{double fibration} in the sense of Gelfand (see \cite{GuSt79}) and Helgason \cite{Helgason}.
In particular,  for each $\vt\in W$ and $x\in X$, 
\be\label{def ztx}
\zvt^x=\big\{y: (x,y)\in\zvt\big\}\subset Y
\ee
is a smooth surface of  codimension $k$.
As in the paper \cite{GGIP12} by Grafakos, Palsson and the first two authors, 
for a compact subset $E\subset\R^d$, we define the (two-point) \emph{$\Phi$-configuration set} of $E$ as
\be\label{def DeltaPhi}
\Delta_\Phi(E)=\big\{\Phi(x,y): x,y\in E\big\}\subset \R^k.
\ee
Thus, if $k=d$ and $\Phi(x,y)=x-y$, then $ \Delta_\Phi(E)=E-E$, while if $k=1$ and  $\Phi(x,y)=|x-y|,\, \Delta_\Phi(E)$ is the distance set of $E$ 
\footnote{In  applications, one can easily localize away from  singularities of $\Phi$ at degenerate configurations (at $x=y$ for this $\Phi$).}.
We give further examples below.
\ss

Our goal is  to find a threshold, $s_0=s_0(\Phi)$, such that if $\hd(E)>s_0$ then $\Delta_\Phi(E)$ has nonempty interior.
Similar questions have been studied in, or are accessible to the methods of, a number of works, e.g., \cite{GIP14,BIT16,CLP16,IT19,K19,IL19}.
(There is of course an extensive literature on the related Falconer  distance problem \cite{Falc86}, and its generalizations to configurations,
where the question is what lower bound on $\hd(E)$ 
ensures that $\Delta_\Phi(E)$ has positive Lebesque measure \cite{W99,Erd06,DGOWWZ18,GIOW18}.)
We will show that a sufficient value of $s_0(\Phi)$  can be expressed in terms of $d,\, k$ and $\az$, the amount 
of smoothing  on
$L^2$-based Sobolev spaces satisfied by  the family of generalized Radon transforms $\rvt$ defined by $\Phi$. 
\ss

There are a number of results in the configuration  literature, such as Henriot-\L aba-Pramanik \cite{HLP15},
Chan-\L aba-Pramanik \cite{CLP16}, and
Fraser-Guo-Pramanik \cite{FGP19}, which make a Fourier dimension assumption on the set $E$.
In some sense, the Sobolev regularity assumption we make here is analogous to a Fourier decay estimate, but on the family of configurations 
rather than on  $E$, since in the translation-invariant case the Sobolev regularity of the $\rvt$ corresponds to the uniform Fourier decay of the 
measures that are the convolution kernels of the $\rvt$. 
This allows our condition on the set $E$ to be solely in terms of the Hausdorff dimension.

\ms

The amount $\az$ of smoothing on $L^2_s(\R^d)$ is defined as follows.
For  each $\vt\in W$, choose a  compactly supported smooth density $d\rho_{\vt}$  
on  $\zvt$ (such as the Leray density induced by $\Phi$),  and define the generalized Radon transform
\be\label{def rt}
\rvt f(x)=\int_{\zvt^x} f(y).
\ee
Thus, the Schwartz kernel of $\rvt$ is $d\rho_{\vt}(x,y)$. It is known that each  $\rvt$ is a Fourier integral operator (FIO), 
of order $-(d-k)/2$ and  associated to a canonical relation $C_{\vt}$, where $C_{\vt}=N^*\zvt' \subset T^*\R^d\times T^*\R^d$,
the conormal bundle of $\zvt$; see Guillemin and Sternberg \cite{GuSt}, and also Phong and Stein \cite{PS86}.

If $C_{\vt_0}$ is a local canonical graph (see the more extended discussion in Sec. \ref{sec back} below), then so are all the $C_{\vt}$
for $\vt$ close to $\vt_0$, and, for $s\in\R$, $\rvt:L^2_s\to L^2_{s+(d-k)/2}$ uniformly as $\vt$ ranges
over   $V$. However, general canonical relations need not be local canonical graphs and there may be a loss of derivatives relative to this estimate;
it is thus useful to describe the amount of smoothing both in terms of the absolute number, $\az$,  of derivatives that the $\rvt$ add on $L^2$, 
and also in terms of the loss, $\beta_\Phi$, relative to the optimal possible smoothing.
Hence, for our first result, we state the regularity  assumption on the $\rvt$ as either of the  equivalent conditions that  
\smallskip

(i) there is an {\it absolute smoothing} $\ap$,  $0\le\az\le (d-k)/2$  such that,  for any $s\in\R$ one has
\be\label{est rt ap}
||\rvt||_{L^2_{s}(\R^d)\to  L^2_{s+\az}(\R^d)} \le C_s,
\ee
with $C_s$ bounded as $\vt$ varies over $W$; or 
\ss

(ii)  for some {\it relative loss} $\bp$,  $0\le \bp \le (d-k)/2$, 
\be\label{est rt bp}
||\rvt||_{L^2_{s}(\R^d)\to  L^2_{s+(d-k)/2-\bp}(\R^d)} \le C_s,
\ee
which holds iff \eqref{est rt ap} does, with $\bp=(d-k)/2-\az$.
Sharp values of $\az$ and $\bp$ are known in a number of degenerate geometries,
but in most for the examples below the canonical relations are nondegenerate, so that  $\bp=0$.
Our first result is the following.
\ss

\begin{theorem} \label{thm same d}
Suppose that $\,\Phi:\R^d\times\R^d\to \R^k$ satisfies \eqref{cond dPhi} 
and that  the generalized Radon transforms $\rvt$ in \eqref{def rt} satisfy \eqref{est rt ap}, \eqref{est rt bp}. 
Then, if $E\subset\R^d$ is a compact set with $\hd(E)>d-\az=\frac{d+k}2+\bp$, it follows that $\Delta_\Phi(E)$ has nonempty interior.
\end{theorem}
\ss

\bc\label{cor MS} 
(Mattila-Sj\"olin \cite{MS99}; Iosevich-Mourgoglou-Taylor \cite{IMT11}) Let $E\subset\R^d$ be a compact set and 
let $\Delta(E)$ be its distance set, $\Delta(E)=\{|x-y|: x,y\in E\}$. More generally, let $\Delta(E)$ be 
defined by the translation-invariant metric induced by a norm $||\cdot ||$  on $\R^d$ 
whose unit sphere has strictly positive Gaussian curvature.
Then, if $\hd(E)>(d+1)/2$, it follows that $\Delta(E)$ contains an open interval.
\ec
\ss

This follows from Thm. \ref{thm same d} by
taking  $k=1$ and  $\Phi(x,y)=||x-y||$, so that $\Delta_\Phi(E)=\Delta(E)$.
For $0<t<\infty$, the generalized Radon transform $\mathcal R_t$ is the spherical mean operator for radius $t$ with respect to $||\cdot||$; as is well-
known (see, for example, \cite{St93}), this is an FIO of order $-(d-1)/2$ associated with a canonical graph, so that, on $L^2_s(\R^d)$,  
$\mathcal R_t$ is smoothing of order $\ap=(d-1)/2$ and $ \bp=0$. 
Hence,  in the same range as in  \cite{MS99}, namely for $s>(d+1)/2$
(also the same range as in Falconer's original result for the distance set problem \cite{Falc86}),
the distance set  has nonempty interior, i.e., contains an open interval.
( For certain specific thin sets below the $(d+1)/2$ threshold, it has been shown that the distance set has nonempty interior; see for example 
\cite{ST17} by Simon and the third author of the current paper.)
\ss

\bc\label{cor riemann}
Let $(M,g)$ be a $d$-dimensional Riemannian manifold, $U\subset M$ an open set 
on which there are no conjugate points, and identify $U$ with an open subset of $\R^d$. Setting $\Phi(x,y)=d_g(x,y)$, 
the Riemannian distance and, for $E\subset U$, the Riemannian distance set $\Delta^g(E):=\{d_g(x,y): x,y \in U\}$,
if $\hd(E)>(d+1)/2$, then $\Delta^g(E)$ contains an open interval.
\ec

This follows from Thm. \ref{thm same d} in the same way as for Cor. \ref{cor MS}:
the operators $\mathcal R_t$ are the Riemannian spherical means, 
which are still FIO of order $-(d-1)/2$  associated with canonical graphs \cite{Ts76}, so that again $\ap=(d-1)/2,\, \bp=0$.
\ms

One can also prove a result for multi-parameter distance sets (cf. \cite{IJP17}): 

\bc\label{cor kfold}
Suppose $d=d_1+\dots+d_k$, with all $d_j>1$, and write $x\in\R^d$ as $x=(x^1,\dots,x^k)$ with $ x^j\in\R^{d_j},\, 1\le j\le k$. 
For $E\subset\R^d$ compact, define
$$\Delta^{(k)}(E)=\left\{\left(|x^1-y^1|,\dots, |x^k-y^k|\right):x,y\in E\right\}\subset\R^k.$$
Then, if $\hd(E)>d-\frac12\min\{d_j-1\}=\frac12\max\{d-d_j+1\}$, it follows that $\Delta^{(k)}(E)$ has nonempty interior in $\R^k$.
\ec

This follows by noting that, for $\vt\in\R^k,\, t_j>0$ for all $j$, 
the operator $\rvt$ is convolution with a product of surface measures on spheres and thus has Fourier multiplier that decays uniformly as $(1+|\xi|)^{-\frac12\min\{d_j-1\}}$ (and no better),
so that Thm. \ref{thm same d} applies with the same $k$ and $\ap=\frac12\min\{d_j-1\}$.
\ms

Thm. \ref{thm same d} is subsumed in a more general version with two sets,
whose proof is essentially the same. 
For $E,\, F\subset \R^d$, define  the {\it$\Phi$-configuration set  of $E$ and $F$} as
\be\label{def asym Delta}
\Delta_\Phi(E,F):=\left\{\Phi(x,y): x\in E,\, y\in F\right\}\subset\R^k.
\ee
Then we have
\begin{theorem} \label{thm EandF}
Suppose that $\,\Phi:\R^d\times\R^d\to \R^k$ satisfies \eqref{cond dPhi} 
and that the generalized Radon transforms $\rvt$ in \eqref{def rt} satisfy \eqref{est rt ap}, \eqref{est rt bp}. 
Then, if $E,\, F\subset\R^d$ are compact sets with $\hd(E)$ and $\hd(F)$ satisfying
 \be\label{cond same d}
\hd(E)+\hd(F)>2d-2\az=d+k+2\bp,
 \ee
it follows that  $\Delta_\Phi(E,F)$ has nonempty interior.
\end{theorem}
The method for proving Thm.  \ref{thm EandF} is flexible and allows one 
to obtain  the nonempty interior of $\Delta_\Phi$ in  {asymmetric} settings, 
where not only are $E$ and $F$ possibly different, but (i) $\Delta_\Phi$ is determined by `points' $x$ and $y$ which may belong to  different spaces, 
of possibly different dimensions; and (ii) $\Phi$ can be manifold-valued.

To make this precise, let $W, X$ and $Y$ be smooth manifolds, of dimensions $k\ge 1$ and $d_1 \ge d_2\ge 1$, resp.; 
note that the notion of Hausdorff dimension for compact subsets of $X$ or $Y$ is well-defined.
Let $\Phi:X\times Y\to W$ be a $C^\infty$ function satisfying the double-fibration condition \eqref{cond dPhi}, 
at least away from a lower dimensional subvariety.
If $E\subset X$, $F\subset Y$ are 
compact,  define
\be\label{def W Delta}
\Delta_\Phi(E,F):=\left\{\Phi(x,y): x\in E,\, y\in F\right\}\subset W.
\ee
Our goal is to find a lower bound on
 $\hd(E)+ \hd(F)$ ensuring that $\Delta_\Phi(E,F)$ 
 has nonempty interior in $W$.
For $\vt\in W$, the operators $\rvt$ are now FIO of order 
$$m:=0+k/2-(d_1+d_2)/4=-(d_2-k)/2-(d_1-d_2)/4,$$
where we have written the order  on the right
(a way that is standard for FIO between spaces of different dimensions)
so as  to the isolate the best possible smoothing order of the operator, namely $(d_2-k)/2$ (see Sec. \ref{sec back} below).
 The analogue of \eqref{est rt ap}, \eqref{est rt bp} is then
\be\label{est rt asymm}
||\rvt||_{L^2_{s}(Y)\to  L^2_{s+\az}(X)}=||\rvt||_{L^2_{s}(Y)\to L^2_{s+\frac{d_2-k}2-\bp}(X)} \le C_s,\,\hbox{ for all }\vt\in W_0\subset W, 
\ee
for some $0\le\az,\, \bp\le (d_2-k)/2$, related by $\ap=(d_2-k)/2-\bp$.
When the $C_{\vt}$ are nondegenerate canonical relations,  $\ap=(d_2-k)/2,\, \bp=0$.
In this  generality,  the following holds, and also implies Thm. \ref{thm EandF} and thus Thm. \ref{thm same d}.

\begin{theorem} \label{thm asymm}
Let $W, X$ and $Y$ be smooth manifolds, of dimensions $k\ge 1$ and \linebreak$d_1 \ge d_2\ge 1$, resp.
Suppose that $\,\Phi:X\times Y\to W$ satisfies  condition \eqref{cond dPhi} 
and  the  $\rvt$ in \eqref{def rt} satisfy \eqref{est rt asymm}  for some $0\le\az,\, \bp\le (d_2-k)/2$.
If $E\subset X$, $F\subset Y$ are compact sets, then $\Delta_\Phi\left(E,F\right)$ defined in \eqref{def W Delta} has nonempty interior in $W$  if
\be\label{cond asymm}
\hd(E)+\hd(F)> d_1+d_2-2\az=d_1+k+2\bp.
\ee
\end{theorem}

\section{Examples}\label{sec examples}

Before proving  Thm. \ref{thm asymm}, we state a series of corollaries  which follow by applying  Thms.   \ref {thm EandF} or \ref{thm asymm} to the relevant $\Phi$, deferring  of some of the details  to Sec. \ref{sec details}.

\subsection{Generalized distance sets}\label{subsec distance}

A diverse class of examples with $k=1$ arises when  the $x$ and $y$ are points, subspaces or even  submanifolds of some ambient $\R^d$
or Riemannian manifold $(M,g)$. 
Suppose $X$ and $Y$ are 
smooth families of closed submanifolds (or possibly subvarieties) 
in $\R^d$ or $(M,g)$, and let $\Phi:X\times Y\to\R$ be the standard distance between subsets,
\be\label{def dist}
\Phi(x,y):=\dist(x,y)=\inf_{a\in x,\, b\in y} |a-b| \quad [\hbox{ or } \inf d_g(a,b)\, ].
\ee
In many cases, $\Phi$ is $C^\infty$ on $X\times Y$, and \eqref{cond dPhi} holds (with $k=1$), away from a lower dimensional variety, corresponding to degenerate configurations.
(Note that, even if $X=Y$, the function $\Phi$ is typically not a metric, since $\Phi(x,y)=0$ does not imply that $x=y$, 
but only that $x\cap y \ne \emptyset$.)
Then, for sets $E\subset X,\, F\subset Y$,  let $\Delta_\Phi(E,F)$ be defined by 
(\ref{def asym Delta}). 
Applications  of  Thm. \ref{thm asymm} to configuration sets of this type include the following: 
\ms

{\bf Lines--points in $\R^2$ and hyperplanes-points in $\R^d$:} Let $Y=\R^2$ and $X=M_{1,2}$, the space of all affine lines in $\R^2$. 
As with the Radon transform \cite{Helgason}, it is convenient to parametrize $M_{1,2}$ as $\sone\times\R$, 
with $(\omega,s)\leftrightarrow \{y\in\R^2: \om\cdot y-s=0\}$. (Since $(-\om,-s)$ and $(\om,s)$ 
correspond to the same line, there is a 2-1 redundancy, but this is harmless.)
Then $\dist((\om,s),y)=|\om\cdot y-s|=:\Phi((\om,s),y)$, which is smooth away from the incidence relation $Z_0$.
The question of interest is:  for collections $E\subset M_{1,2}$ and $F\subset\R^2$  of lines and points, resp., 
what lower bounds on  $\hd(E)$ and $\hd(F)$ ensure that 
$\Delta_\Phi(E,F)$ contains an interval? 
(Finite field Falconer-type problems for this geometry  were studied in \cite{V14,PPSVV17,BIP17,PPV19,B19}.)

The operators $\mathcal R_t$ are sums of translates by $\pm t$ of the Radon transform  in the $s$ variable
so that the $C_t$ are local canonical graphs and $\bp=0$.
Similarly, for $d\ge 3$, let $Y=\R^d$ and $X=M_{d-1,d}$, the Grassmannian of all affine hyperplanes in $\R^d$. 
As for $d=2$, we may parametrize $M_{d-1,d}$ as $\sd\times\R$
and, with $k=1$, take $\Phi((\om,s),y):=\dist((\om,s),y)=|\om\cdot y-s|$.
Then the $C_t$ are again local canonical graphs, so that $\bp=0$ and  a consequence of Thm. \ref{thm asymm} is
\bc\label{cor hyperplanespoints}
If $E\subset M_{d-1,d},\, F\subset\R^d$ with $\hd\left(E\right)+\hd\left(F\right)>d+1$,
then $\Delta_\Phi(E,F)$ contains an interval.
\ec
\ms

{\bf Circles--points in $\R^2$ and spheres-points in $\R^d$:} For  a family of non-equidimen-sional examples, again let $Y=\R^2$ and now let $X=S_2$ be the  space of all circles in $\R^2$,
which is 3-dimensional. We may parametrize $S_2$  by the center and radius of each circle, $S_2=\{(a,r):a\in\R^2,\, r>0\}$.
Then $\Phi((a,r),y):=|y-a|-r$ defines the circle-point relation; for each $t\in\R$, $Z_t^{\pm}=\{((a,r),y): |y-a|=r\pm t\}$
is smooth away from the singular set at $r=\mp t$. The operator $\rt\in I^{-\frac12-\frac14}(C_t)$ where $C_t$, 
which is the same canonical relation as for the forward solution of the wave equation (with $r$ playing the role of time), 
but translated by $\pm t$ in the $r$ variable.
Since this is nondegenerate, $\bp=0$
and  so it follows that if $E\subset S_2$, a set of circles, and $ F\subset\R^2$, a set of points,  satisfy
$\hd(E)+\hd(F)>4$, then Thm. \ref{thm asymm} applies. This extends to higher dimensions:

\bc\label{cor circlespoints}
If $S_d$ is the  $(d+1)$-dimensional space of  all spheres in $\R^d$, and $E\subset S_d,\, F\subset\R^d$ satisfy
$\hd(E)+\hd(F)>d+2$, then
$\Delta_\Phi(E,F)\subset\R$ has nonempty interior. 
\ec
\ms

{\bf Lines-points in $\R^d,\, d\ge 3$:} Let $X=M_{1,d}$ be the $(2d-2)$-dimensional Grassmannian of all affine lines in $Y=\R^d$. 
For $L\in M_{1,d}$ and $y\in\R^3$, let $\Phi(L,y)=\dist(L,y)$. 
One can show (see Sec. \ref{sec details}) that for $t>0$, $\rt$ is an FIO with a nondegenerate canonical relation and thus $\bp=0$. 
From Thm. \ref{thm asymm} we obtain

\bc\label{cor linespoints}
If  $E\subset M_{1,d}$ and $F\subset\R^d$ satisfy $\hd(E)+\hd(F)>2d-1$, 
then  $\Delta_\Phi(E,F)$ has nonempty interior.
\ec
\ms

{\bf Lines-lines in $\R^d,\, d\ge 3$:} 
Now let $X=Y=M_{1,d}$. For $L,\, L'\in M_{1,d}$, the distance  $\dist(L,L')$
as  in (\ref{def dist}) is positive if $L\cap L'=\emptyset$, and is a smooth function, satisfying the double fibration condition, away 
from the lower dimensional incidence variety  $\{(L,L'):\, L\cap L'\ne\emptyset\}$. 
For $t>0$, the $\rt$  are FIO associated with  canonical graphs, so that $\bp=0$, yielding
\bc\label{cor lineslines}
If $E,\, F\subset M_{1,d}$ are sets of lines with $\hd(E)+\hd(F)>2d-1$, then
$\Delta_\Phi(E,F)$ has nonempty interior.
\ec

\subsection{Higher dimensional configuration sets}\label{subsec higher} 
Values of  $k$ greater than 1 arise  when configurations are encoded by vector-valued
(or manifold-valued) data.
\ms

{\bf Configurations determined by ensembles of quadratic forms:} Let $Q_1,\dots,Q_k$ be quadratic forms on $\R^d$. 
Define the ensemble $\vec{Q}=(Q_1,\dots,Q_k)$ to be {\it nonsingular} if
$c_1Q_1+\cdots +c_kQ_k$ is nonsingular
for all $\vec{c}=(c_1,\dots,c_k)\in\R^k\setminus {\mathbf 0}$. 
The sharp restrictions on $k$ and $d$ in order  that such ensembles to exist were found by Adams, Lax and Phillips \cite{ALP65}: 
there exists such a $k$-dimensional nonsingular family of quadratic forms on $\R^d$ iff $k\le \rho(d/2)+1$,
where $\rho(\cdot)$ are  {Radon-Hurwitz numbers}, defined by
\ms

\ni {\bf Def.}  
If $n$ is a positive integer, write $n=(2l+1)\cdot 2^m$, the factorization of $n$ into the product of an odd integer and an integral power of 2. Express $m$ modulo 4 as $m=p+4q,\, 0\le p\le 3$. Then the {\it Radon-Hurwitz number} of $n$ is $\rho(n):=2^p+8q$.
If $n$ is a half-integer, then  $\rho(n):=0$. 
\smallskip

\bc\label{cor quadratic}
If $\vec{Q}$ is  a nonsingular family on $\R^d$, define $\Phi:\R^d\times \R^d\to \R^k$ by $\Phi(x,y)= \vec{Q}(x-y)$. 
Then if $E,\, F\subset\R^d$ are compact and
$\hd(E)+\hd(F)>d+k$, it follows that
$$\Delta_{\Phi}(E,F)=\big\{\left(Q_1\left(x-y\right),\dots,Q_k\left(x-y\right)\right)\in\R^k: x\in E,\, y\in F\big\}$$
has nonempty interior in $\R^k$.
\ec

If $d$ is odd, only $k=1$ is possible. 
For $d=2$, $\rho(d/2)+1=2$, and $Q_1(x)=x_1^2-x_2^2,\, Q_2(x)=x_1x_2$ is a nonsingular ensemble,
but  for this $\vec{Q}$   Thm. \ref{thm EandF} is vacuous, since it requires $\hd(E)+\hd(F)>4$. 
However, note that for $d=4$, one has $\rho(d/2)+1=2+1=3$, so that there exist nonsingular ensembles  for which Cor. \ref{cor quadratic} yields a nontrivial result; 
e.g., with $k=3$, if $\hd(E)>7/2$ then $\Delta_\Phi(E,E)$ has nonempty interior.
\ms

{\bf Heisenberg circles and spheres:}
Let $\H^1$ be the 3-dimensional  Heisenberg group, $\H^1\simeq \R^2\times\R$ with product
$$(x',x_3)\cdot (y',y_3)=\big(x'+y', x_3+y_3+\frac12(x')^TJy'\big),\quad
J=\left[\begin{matrix} 0 & -1 \cr 1 & 0 \end{matrix}\right].$$
Define $\Phi:\H^1\times\H^1\to \R^2$ by
\be
\label{eqn orig phi for heis}
\Phi(x,y)=\vt= \left(|x'-y'|, x_3-y_3+\frac12\left(x_1y_2-x_2y_1\right)\right) \in\R^2, 
\ee
where $t_1$ is the distance from the origin in the first two coordinates of $x\cdot y^{-1}$ and $t_2$ is the `height' of 
$x\cdot y^{-1}$ above the  `plane of good directions', $\{x_3=0\}$. 
For each $\vt\in \R_+\times\R$, the corresponding $\rvt$ is the generalized Radon transform which averages over 
group translates of the origin-centered circle of radius ${t_1}\subset\{x_3=0\}$, translated in the central direction by $t_2$.
Such operators (for $t_2=0$) have been studied by Nevo and  Thangavelu \cite{NT97} and M\"uller and Seeger \cite{MuSe04}.
Here, in the notation of Thm. \ref{thm EandF}, $d=3,\, k=2$, and  it is known (see \cite{MuSe04} 
and Sec. \ref{sec details}) that each $C_{\vt}$ is a two-sided fold 
(or folding canonical relation in the sense of Melrose and Taylor \cite{MeTa85}), so that $\bp=1/6$ by \cite{MeTa85}. 
Hence, from Thm. \ref{thm EandF} we obtain

\bc\label{cor heis}
If $E,\, F\subset\H^1$ with 
$\hd(E)+\hd(F)> 16/3$, then $\Delta_\Phi(E,F)\subset\R^2$ has nonempty interior.
\ec
\ss

This  generalizes to the more general setting of \cite{MuSe04}, with $\H^1$ replaced by a nondegenerate step-two nilpotent group $G$
and $\mathbb S^1$ replaced by a  hypersurface $\Sigma$ with strictly positive Gaussian curvature in the step 1 space: 
Suppose $\g=\g_1\oplus \g_2$ is a graded nilpotent Lie algebra, with $dim(\g_1)=n,\, dim(\g_2)=m$, and $G=exp(\g)\simeq \R^{n+m}$ is  the 
resulting connected, simply-connected Lie group, with projections $\pi_1:G\to \g_1\simeq\R^n$ and $\pi_2:G\to \g_2\simeq\R^m$.
We assume that $\g$ satisfies the nondegeneracy condition that the skew-linear form
\be\label{def nondeg}
B_\zeta(u,v):=\langle \zeta,\, [u,v] \rangle \hbox{ is nondegenerate on } \g_1,\quad\forall \zeta\in \g_2^*\setminus 0.
\ee
(In particular, $n$ must be even.) 
 Let $||\cdot ||$ be a $C^\infty$ norm on $\g_1$ whose unit sphere  $\Sigma$ has everywhere positive Gaussian curvature,
and  define $\Phi:G\times G\to  \R^{m+1}$ by
\be\label{def phi nilpotent}
\Phi(x,y)=\vt:=\big(\, ||\pi_1(x\cdot y^{-1})||, \pi_2 (x\cdot y^{-1})\, \big).
\ee
Then, one can show that the the canonical relations $C_{\vt}$ of the $\rvt$ are associated with two sided folds, so that again $\bp=1/6$,
and the required lower bound on $\hd(E)+\hd(F)$ is $(n+m)+(m+1)+\frac13=n+2m+\frac43$,
implying

\bc\label{cor step two}
If $E,\, F\subset G$ with $\hd(E)+\hd(F)>n+2m+\frac43$,
then $\Delta_\Phi(E,F)$ has nonempty interior.
\ec

\ms

{\bf Configurations determined by a curve:} Let $\gamma:\R\to\R^d$ be a smooth curve. Writing $x=(x_1,x')$, etc., 
we assume $\gamma$  is of the form $\gamma(\tau)=(\tau,g(\tau)),\, g=(g_2,\dots,g_d):\R\to\R^{d-1}$. Define $\Phi:\R^d\times\R^d\to\R^{d-1}$ by 
\be\label{eqn moment}
\Phi(x,y)=x'-y'-g(x_1-y_1),
\ee
which measures the displacement of $x$ off of $y+\gamma$. 
All of the operators $\rvt$ are translates of $\mathcal R_{\bf 0}$ and thus have the same structure and satisfy the same estimates. For example, if $\gamma$ is the moment curve $\gamma(\tau)=(\tau,\tau^2,\dots,\tau^d)$, then $\ap=\frac1d$, and from Thm. \ref{thm EandF} we obtain

\bc\label{cor curves}
If $\gamma$ is the moment curve and compact sets $E,\, F\subset \R^d$ satisfy $\hd(E)+\hd(F)>2d-\frac{2}d$,
then $\Delta_\Phi(E,F)$ has nonempty interior.
\ec

\section{Background}\label{sec back}

We give a brief survey of the Fourier integral operator theory needed, referring to  H\"ormander \cite{Hor71,Hor85} for more background and further details. 

Let $X$  be a smooth manifold of dimension \nolinebreak$d$.
The cotangent bundle $T^*X$ is a symplectic manifold with respect to the {\it canonical two-form}, $\omega=\sum d\xi_j\wedge dx_j$ 
(with respect to any local coordinates). We denote  the {\it zero-section} of $T^*X$, $\{\xi=0\}$, by $\z$.
A {\it conic Lagrangian submanifold} of $T^*X$ is  a smooth, conic (i.e., invariant under $(x,\xi)\to (x,\tau\xi)$ for $0<\tau<\infty$) submanifold $\Lambda\subset T^*X\mzs$ of dimension $\dim(\Lambda)=d=\frac12 \dim(T^*X)$ such that $\omega|_\Lambda\equiv 0$.
\ms

Now let $X$ and $Y$ be smooth manifolds of dimensions $d_1,\, d_2$, resp. Then $T^*X,\, T^*Y$ are each symplectic manifolds, with canonical two-forms we denote by $\omega_{T^*X},\, \omega_{T^*Y}$, resp. Equip $T^*X\times T^*Y$ with the {\it difference symplectic form}, $\omega_{T^*X} - \omega_{T^*Y}$. For our purposes, a {\it canonical relation} will mean a  submanifold, 
$C\subset (T^*X\mzs)\times(T^*Y\mzs)$  (hence of dimension $d_1+d_2$), which is  conic Lagrangian with respect to $\omega_{T^*X} - \omega_{T^*Y}$.
\ms

For some $N\ge 1$, let $\phi: X\times Y\times (\R^N\mzs)\to \R$ be a smooth phase function 
which is positively homogeneous of degree 1 in $\theta\in\R^N$, i.e., 
$\phi(x,y,\tau\theta)=\tau\cdot\phi(x,y,\theta)$ for all $\tau\in\R_+$. 
Let $\Sigma_\phi$ be the {\it critical set} of $\phi$ in the $\theta$ variables,
$$\Sigma_\phi:=\{(x,y,\theta)\in X\times Y\times (\R^N\mzs): d_\theta\phi(x,y,\theta)=0\},$$
and 
$$C_\phi:=\{(x,d_x\phi(x,y,\theta); y, -d_y\phi(x,y,\theta)): (x,y,\theta)\in\Sigma_\phi\},$$
both of which are conic sets.
If we impose  the first order nondegeneracy conditions
$$d_{x}\phi(x,y,\theta)\ne 0\hbox{ and } d_{y}\phi(x,y,\theta)\ne 0, \forall (x,y,\theta)\in \Sigma_\phi,$$
then $C_\phi\subset (T^*X\mzs)\times(T^*Y\mzs)$.
If in addition one demands that
$$\rank[d_{x,y,\theta}d_\theta\phi(x,y,\theta)]=N,\, \forall\, (x,y,\theta)\in\Sigma_\phi,$$
then $\Sigma_\phi$ is smooth, $dim(\Sigma_\phi)=d_1+d_2$, and the map
\be\label{eqn cphi param}
\Sigma_\phi\owns (x,y,\theta)\to \left(x,d_x\phi\left(x,y,\theta\right); y,-\d_y\phi\left(x,y,\theta\right)\right)\in C_\phi
\ee
is an immersion, whose image is an immersed canonical relation; the phase function $\phi$ is said to {\it parametrize} $C_\phi$.
\ss

For a canonical relation $C\subset (T^*X\mzs)\times(T^*Y\mzs)$  and $m\in\R$,
one  defines $I^m(X,Y;C)=I^m(C)$, the class of  {\it Fourier integral operators} $A:\mathcal E'(Y)\to \mathcal D'(X)$ of order \nolinebreak$m$,
as the collection of operators whose Schwartz kernels are locally finite sums of oscillatory integrals of the form
$$K(x,y)=\int_{\R^N} e^{i\phi(x,y,\theta)} a(x,y,\theta)\, d\theta,$$
where $a(x,y,\theta)$ is a symbol of order $m-\frac{N}2+\frac{d_1+d_2}4$ and $\phi$ is a phase function as above, 
parametrizing some $C_\phi\subset C$.
\ms

The  FIO relevant for this paper are the {\it generalized Radon transforms} $\rvt$ 
determined by  defining functions $\Phi:X\times Y\to \R^k$ 
satisfying \eqref{cond dPhi}.   The Schwartz kernel of each $\rvt$ is a smooth multiple of $\delta_k(\Phi(x,y)-\vt)$, 
where $\delta_k$ is the delta distribution on $\R^k$. From the Fourier inversion representation of $\delta_k$, we see that $\rvt$ has kernel
$$K_{\vt}(x,y)=\int_{\R^k} e^{i(\Phi(x,y)-\vt)\cdot\theta} \, b(x,y)\cdot 1(\theta)\, d\theta,$$
where $b\in C_0^\infty$. Since the amplitude is a symbol of order 0, $\rvt$ is an FIO of order $0+\frac{k}2-\frac{d_1+d_2}4=-\frac{d_1+d_2-2k}4$
associated with the canonical relation parametrized as in \eqref{eqn cphi param}
by $\phi(x,y,\theta)=(\Phi(x,y)-\vt)\cdot\theta$, which is the twisted conormal bundle
of the incidence relation $\zvt$,
$$C_{\vt}=N^*\zvt':=\big\{\big(x,\sum_{j=1}^k d_x\Phi_j(x,y)\theta_j; 
y,-\sum_{j=1}^k d_x\Phi_j(x,y)\theta_j\big): (x,y)\in \zvt,\, \theta\in\R^k\mzs\big\}.$$
For $W$-valued defining functions $\Phi$,
 as in Thm. \ref{thm asymm}, this discussion is modified slightly by introducing local coordinates on $W$.
\ss

For a general canonical relation, $C$, the natural projections $\pi_L:T^*X\times T^*Y \to T^*X$ and $\pi_R:T^*X\times T^*Y \to T^*Y$ restrict to 
$C$, and by abuse of notation we refer to the restricted maps with the same notation.
One can show that, at any point $c_0=(x_0,\xi_0;y_o,\eta_0)\in C$, one has $\corank(D\pi_L)(c_0)=\corank(D\pi_R)(c_0)$; we say that
the canonical relation $C$ is {\it nondegenerate} if this corank is zero at all points of $C$, i.e., if
$D\pi_L$ and $D\pi_R$ are of maximal rank. 
If $\dim(X)=\dim(Y)$, then $C$ is nondegenerate iff $\pi_L,\, \pi_R$ are local diffeomorphisms, and then $C$ is a {\it local canonical graph}, i.e., locally near any $c_0\in C$ equal to the graph of a canonical transformation. If $\dim(X)=d_1>d_2=\dim(Y)$, then $C$ is nondegenerate iff $\pi_L$ is an immersion and $\pi_R$ is a submersion.
To describe the $L^2$-Sobolev estimates for FIO associated with $C$, it is convenient to normalize the order by 
considering $A\in I^{m-\frac{|d_1-d_2|}4}$. One has

\begin{theorem}\label{thm nondeg fio}
Suppose that $\dim(X)=d_1,\, \dim(Y)=d_2$, $C\subset  (T^*X\mzs)\times(T^*Y\mzs)$ is a nondegenerate canonical relation, 
and $A\in I^{m-\frac{|d_1-d_2|}4}$ has a compactly supported Schwartz kernel. Then $A:L^2_s(Y)\to L^2_{s-m}(X)$ for all $s\in\R$.
\end{theorem}

In particular, if $C$ is a local canonical graph, then $A\in I^m\implies A:L^2_s\to L^2_{s-m}$; this is relevant to a number of the Corollaries above.

On the other hand, for Cors. \ref{cor circlespoints} and \ref{cor linespoints}, 
 the canonical relations cannot be canonical graphs,
 since the dimensions of $X$ and $Y$ differ, but the canonical relations are nondegenerate.
For Cor. \ref{cor circlespoints},  parametrized by
the pair of phase functions $\phi^\pm((a,r),y;\theta)=(|y-a|-r\mp t)\theta$ on $\R^{d+1}\times\R^d\times (\R\mzs)$,
$C_t$ is, away from $r\mp t =0$, given by
$$C_t=\{(a,r,\theta\omega,\mp|\theta|;a+(r\pm t)\omega,\theta\omega): (a,r)\in\R^{n+1},\omega\in\mathbb S^{d-1},\theta\in\R\mzs\}.$$
One sees by inspection that $D\pi_R$ has rank $2d$ everywhere and thus $C_t$ is nondegenerate away from $\{r\mp t=0\}$ (in fact, its natural 
extension across those points is also  smooth and nondegenerate, but for our purposes we will not need that).
Thus, when localized away from $r=\pm t$, the operators $\rt$ are in $I^{-\frac12-\frac14}(C_t)$ 
and by Thm. \ref{thm nondeg fio} map $L^2_s(\R^d)\to L^2_{s+(1/2)}(\R^{d+1})$.
\ss

Returning to general $C$, if the corank of $D\pi_L$ (and thus that of $D\pi_R$) is $\le k$ at all points of $C$, then  an $A\in I^m(C)$ maps 
no worse than $L^2_s(Y)\to L^2_{s-m-(k/2)}(X)$.
However, for classes of $C$ for which $\pi_L$ and/or $\pi_R$ degenerate in specific ways, 
the loss of derivatives is often less than $k/2$. The first and best known result of this type is the following, which
we use in the analysis below of Cors. \ref{cor heis}  and \ref{cor step two}.
Suppose $\dim(X)=\dim(Y)=d$ and at any degenerate points  $c_0\in C\subset  (T^*X\mzs)\times(T^*Y\mzs)$, both $\pi_L$ and $\pi_R$ have 
Whitney fold singularities. Such  $C$ were introduced by Melrose and Taylor \cite{MeTa85} and called {\it folding canonical relations}
(also called {\it two-sided folds} \cite{GS94}). 

\begin{theorem}\label{thm mt}
\cite{MeTa85} 
If $C$ is a folding canonical relation and $A\in I^m(X,Y;C)$ 
with compactly supported Schwartz kernel,
then $A:L^2_s(Y)\to L^2_{s-m-(1/6)}(X),\, \forall s\in\R$.
\end{theorem}

\section{Proof of  theorem  \ref{thm asymm}}\label{sec proofs}

To begin the proof, recall\footnote{We refer to the monographs of Mattila \cite{Mat95,Mat15} for background definitions and results.} 
that if $E\subset\R^d$ is a compact set with $\hd(E)>s$, then
there exists a {Frostman measure} on $E$: 
a probability measure $\mu$, supported on $E$ and of finite $s$-energy:

$$\int_E \int_E |x-y|^{-s}\, d\mu(x)\, d\mu(y) <\infty,$$
or equivalently, 
\be\label{eqn energy}
\int_E |\hat{\mu}(\xi)|^2\cdot |\xi|^{s-d}\, d\xi<\infty.
\ee
Since $\mu\in\mathcal E'(\R^d)$, $\hat{\mu}\in C^\omega$ and thus it follows from \eqref{eqn energy} that $\mu\in L^2_{(s-d)/2}(\R^d)$.
This last fact also holds in the more general setting of $E\subset X$,  a compact subset of a $d$-dimensional manifold $X$ with $\hd(E)>s$.
\ms

Now, in the context of Thm. \ref{thm asymm}, suppose that $dim(X)=d_1,\, dim(Y)=d_2$,  $E\subset X,\, F\subset Y$, 
with  $\hd(E)+\hd(F)>d_1+k+2\bp$. Then we can find $s_1,\, s_2$ such that
$\hd(E)>s_1,\, \hd(F)>s_2$, also 
satisfying 
\be\label{eqn s cond}
s_1+s_2>d_1+k+2\bp
\ee
If $\Phi:X\times Y\to W$ satisfies \eqref{cond dPhi}, then
choices of Frostman measures $\mu_1,\, \mu_2$  on $E,\, F$, relative to $s_1,\, s_2$, resp.,
induce a  \emph{configuration measure}  $\nu$ on  $W$, defined by
$$\nu(A)=\int_E \mu_2\big(\{y\in F\, :\, \Phi(x,y)\in A\}\big)\, d\mu_1(x),$$
or equivalently, for $g\in C_0(W)$,

$$\int_{W} g(\vt) \, d\nu = \int\int_{E\times F} g(\Phi(x,y))\, d\mu_1(x)\, d\mu_2(y).$$

We claim that $\nu$ is absolutely continuous (with respect to Lebesgue measure in any coordinate system on $W$), 
and its density   is a continuous function, $\nu(\vt)$.
First, note that one can formally write
\be\label{eqn nu(t)}
\nu(\vt)=\langle\,\rvt\, \mu_2,\mu_1\,\rangle
\ee
with the pairing on the right hand between elements of Sobolev spaces on $X$. 
Since $\mu_1\in L^2_{(s_1-d_1)/2}(X)$ and $\mu_2\in L^2_{(s_2-d_2)/2}(Y)$, 
and the hypothesis \eqref{est rt asymm} for Thm. \ref{thm asymm} is that, for any $s\in\R$,  $ \rvt:{L^2_{s}(Y)\to  L^2_{s+(d_2-k)/2-\bp}(X)}$ uniformly in $\vt$, the sum of the Sobolev orders of the left- and right-side terms in \eqref{eqn nu(t)} is
\be\label{eqn s for proof}
\frac{s_1-d_1}2+\frac{s_2-d_2}2+\frac{d_2-k}2-\bp =\frac12(s_1+s_2-d_1-k-2\bp)>0,
\ee
with the inequality due to \eqref{eqn s cond}. Thus, the integral representing \eqref{eqn nu(t)} in 
terms of $\widehat{\mu_1},\, \widehat{\rvt\mu_2}$ 
is absolutely convergent by Cauchy-Schwarz, and  by continuity of the integral it depends continuously on the parameter $\vt$.
\ss

To make this rigorous we argue as follows, restricting the analysis to the case when $\Phi:X\times Y\to\R^k$; 
the proof  extends to general $W$ using local coordinates on $W$. 
For a $\chi\in C_0^\infty(\R^k)$  supported in a sufficiently small ball, $\chi\equiv 1$ near $\mathbf 0$,  and with $\int\chi\, d\vt=1$,
set $\chi^\epsilon(\vt):=\epsilon^{-k}\chi(\frac{\vt}{\epsilon})$ the associated approximation to the identity, 
which converges to $ \delta(\vt)$ weakly as $\epsilon\to 0^+$. 
Define $\rvt^\epsilon$ to be the operator with Schwartz kernel  $K_{\vt}^\epsilon(x,y):= \chi_\epsilon(\Phi(x,y)-\vt)$.
Then $\rvt^\epsilon\, \mu_2\in C^\infty(X)$ and depends smoothly on $\vt$, and
thus we can represent $\nu$ as  the weak limit of absolutely continuous measures with smooth densities:
$$\nu=\lim_{\epsilon\to 0^+} \nu^\epsilon:=\lim_{\epsilon\to 0^+}\langle\,\rvt^\epsilon\, \mu_2,\mu_1\,\rangle.$$
Now, the operators $\rvt^\epsilon\in I^{-\infty}(C_{\vt})$, with symbols which converge in the Fr\'echet topology on the space of symbols  
as $\epsilon\to 0$ to the symbol of $\rvt$.
Since the singular limits $\rvt$ satisfy \eqref{est rt asymm}, so do the $\rvt^\epsilon$ uniformly in $\epsilon$. 
Hence, $\nu(\vt)$, being the uniform 
limit of smooth functions of $\vt$, is continuous. 
Furthermore, since $\epsilon^k\cdot \chi^\epsilon$ is bounded below by a constant times the 
characteristic function of the ball of radius $\epsilon$ in $\R^k$,   with constant $C_\Phi$ uniform in $\vt$ we have that
\be\label{eqn epsilon est}
\nu\left(B\left(\vt,\epsilon\right)\right):=
(\mu_1\times\mu_2)\left(\left\{\left(x,y\right): \left|\Phi(x,y)-\vt\right|<\epsilon\right\}\right)\le C_\Phi\epsilon^k.
\ee

So far, we have shown that $\nu(\vt)$ is continuous, so that it is positive on an open set.
Hence, $\Delta_\Phi(E,F)$ is open; to conclude the proof, we  need to show that it is nonempty.
However, this follows because,  as a further consequence of the analysis above, 
it follows  that $\Delta_\Phi(E,F)$ has positive $k$-dimensional Lebesque measure.
In fact, since if $\big\{B(\vt_j,\epsilon_j)\big\}$ is any cover of $\Delta_\Phi(E,F)$, we have
$$1=\mu_1(E)\cdot\mu_2(F)=(\mu_1\times\mu_2)(E\times F)\le\sum_j (\mu_1\times\mu_2)(\Phi^{-1}(B(\vt_j,\epsilon_j))$$
$$=\sum_j \nu\left(B\left(\vt_j,\epsilon_j\right)\right)
\le C_\Phi\sum_j\epsilon_j^k$$
by \eqref{eqn epsilon est}, so that $\sum_j |B(\vt_j,\epsilon_j)|_k\ge C_\Phi'$ is bounded below.
Hence $\Delta_\Phi(E,F)$  has positive $k$-dimensional Lebesgue measure and is therefore nonempty; by the first part of the proof, it in fact has nonempty interior. Q.E.D.

\section{Details of the examples and corollaries}\label{sec details}

This section contains  calculations and additional details to show how some of the 
Corollaries in Sec. \ref{sec examples} follow from the Theorems.

\subsection{Mattila-Sj\"olin and generalizations (Cors. \ref{cor MS} and \ref{cor riemann})}

The  results of Mattila-Sj\"olin \cite{MS99} and Iosevich-Mourgoglou-Taylor \cite{IMT11} follow immediately
from Thm. \ref{thm same d} as indicated in the discussion below Cor. \ref{cor MS}: 
For $t>0$, the spherical mean operators $\rt$ are FIO of order $-(d-1)/2$ associated to canonical relations which  are (under the various 
assumptions) local canonical graphs and thus map $L^2_s\to L^2_{s+(d-1)/2}$; 
furthermore, by standard facts about the dependence on symbols and canonical relations, the operator norms (for fixed $s$) are uniform as $t$ 
ranges over any compact subinterval of $(0,\infty)$. For the Riemannian setting of Cor. \ref{cor riemann}, 
one uses the fact that the same results hold within the conjugate locus \cite{Ts76}.
\ms

Similarly, as indicated in the discussion above its statement, the analysis behind Cor. \ref{cor hyperplanespoints} concerning distances from points to 
lines in $\R^2$ or hyperplanes in $\R^3$ consists of standard facts about the $L^2$ regularity of the Radon transform.

\subsection{Spheres-points (Cor. \ref{cor circlespoints}).} Let $S_d=\{(a,r): a\in\R^d,\, r>0\}\subset \R^{d+1}$ denote
the $(d+1)$-dimensional space of spheres in $\R^d$, and
$Z_t^{\pm}$ be  the smooth points of $\{((a,r),y)\in S_d\times \R^d: |y-a|=r\pm t\}$. 
For $t=0$, the twisted conormal bundle 
$$C_0=N^*Z_0'\subset (T^*\R^{d+1}\setminus\mathbf 0)\times(T^*\R^d\setminus\mathbf 0)$$ 
is  the same canonical relation as for the  solution operator mapping the Cauchy data on $\R^d$  to the solution of the wave equation on $\R^{d+1}$,
which is well-known  to be nondegenerate (see, e.g., \cite{MSS92,Sogge}). 
For $t>0$, $C_t$ is the union of two copies of $C_0$, translated by $\mp t$ in the $r$ variable, and thus is also nondegenerate.
By Thm. \ref{thm nondeg fio}, if $A\in I^{m-\frac14}(C_t)$,  then $A:L^2_s(\R^d)\to L^2_{s-m}(\R^{d+1})$, with norm that is bounded
above  as $t$ ranges over a compact interval in $(0,\infty)$. Hence, \eqref{est rt asymm} is satisfied with $\bp=0$,
and so  Thm. \ref{thm asymm} applies with $d_1=d+1,\, d_2=d,\, k=1$ and $\bp=0$;
thus, if $E\subset S_d$  and $F\subset\R^d$ with $\hd(E)+\hd(F)>(d+1)+1+2\cdot 0=d+2$, then $\Delta_\Phi(E,F)$ has nonempty interior.

\subsection{Lines-points in $\R^d$ (Cor. \ref{cor linespoints})} 
We identify $M_{1,d}$, the Grassmannian of oriented affine lines in $\R^d$, with the tangent bundle $\tsd$:  $(\omega, v)\in \tsd$ 
corresponds to the line $L_{\omega,v}:=\{v+s\omega: s\in\R\}$. 
(Here, we identify $v\in T_\omega\sd$ with the vector $v\in\R^d,\, v\perp \omega$.)
Then define 
$$\Phi((\omega,v),y)=\frac12 dist(y,L_{\omega,v})^2=\frac12 dist(y-v,v^\perp)^2$$
$$= \frac12\big(|y-v|^2-\left(\left(\left(y-v\right)\cdot\omega\right)\omega\right)^2\big)
=\frac12\big|\Pi_\omega^{\perp}(v-y)\big|^2,$$
where $\Pi_\omega^{\perp}$ denotes orthogonal projection onto $\omega^{\perp}$. Then
$$d\Phi'((\omega,v),y)=\big(-((y-v)\cdot\omega)i^*_\omega(y-v),\Pi_\omega^\perp(v-y);\, -\Pi_\omega^{\perp}(v-y)\big),$$
where $i_\omega$ denotes the inclusion of $T_\omega\sd\hookleftarrow \R^d$ and $i^*_\omega$ its transpose.
Note that \linebreak$\Pi_\omega^\perp(v-y)\ne 0$ iff $y\notin\R\cdot\omega$.
Thus, for $t>0$, 
$$Z_t=\{((\omega,v),y)\in M_{1,d}\times \R^d: dist(y,L_{\omega,v})=t\}=\{\Phi=t^2\}$$ 
is smooth and satisfies \eqref{cond dPhi} for 
$$((\omega,v),y)\notin Z_t^{sing}:=\{((\omega,v),y): y\in L_{\omega,v}\},$$
since $\Pi_\omega^\perp(v-y)\ne 0$  at those points.
As coordinates on the $(3d-3)$-dimensional $Z_t\setminus Z_t^{sing}$ 
we may use $y\in\R^d, \omega\in\sd$ and $v\in\{v\in T_\omega\sd: |v|=t\}$.
Letting $\theta\in\R\setminus 0$ be the additional cotangent variable on $C_t=N^*(Z_t\setminus Z_t^{sing})'$, 
the projection $\pi_R:C_t\to T^*\R^d$ is
$$(y,\omega,v,\theta) \to (y, -\theta\cdot\Pi_\omega^\perp(v-y)),$$ 
so that
$$rk(D\pi_R)=(d+1)+rk\Big(\frac{D(\Pi_\omega^\perp(v-y))}{D(\omega,v)}\Big)=(d+1)+(d-1)=2d.$$
Thus, $C_t$ is nondegenerate and Thm. \ref{thm asymm}, with $d_1=2d-2,\, d_2=d,\, k=1$ and $\bp=0$, 
implies that if $E\subset M_{1,3}$ and $F\subset\R^3$ are compact with
$$\hd(E)+\hd(F)>(2d-2)+1+0=2d-1,$$
then $\Delta_\Phi(E,F)$ has nonempty interior.

\subsection{Lines-lines (Cor. \ref{cor lineslines})} Using the same notation as for lines-points, 
we now consider pairs of lines in $\R^d$, say $L=L_{\omega,v},\, L':=L_{\omega',v'}$ with $(\omega,v),\, (\omega',v')\in \tsd$. 
Then, parametrizing $L'$ by $s'\to v'+s'\omega'$, we have
$$dist(L,L')^2=\inf_{y'\in L'} \, dist(y',L)^2=  \inf_{y'\in L'} \, \left|\Pi_\omega^\perp\left(v-\left(v'+s'\omega'\right)\right)\right|^2.$$
A calculation shows that the critical point of the quadratic function of $s'$ in the last expression is $s'=(v-v')\cdot\omega'$, and thus
$$\Phi((\omega,v),(\omega',v')):=\frac12 dist(L,L')^2=\frac12  \big|\Pi_\omega^\perp \Pi_{\omega'}^\perp (v-v')\big|^2$$
$$=\big|v-v'-((v-v')\cdot\omega')\omega' -((v-v')\cdot\omega)\omega+((v-v')\cdot\omega')(\omega\cdot\omega')\omega\big|^2.$$
By a slightly more complicated calculation than in the lines-points case, one sees that this satisfies the double fibration condition away from a lower-
dimensional singular set, and the canonical relations $C_t$ are nondegenerate. Hence by Thm. \ref{thm asymm}, again
with $d_1=2d-2,\, k=1,\, \bp=0$, it follows that if $E,\, F\subset M_{1,d}$ are compact sets of lines with $\hd(E)+\hd(F)>2d-1$, 
then $\Delta_\Phi(E,F)$ has nonempty interior.

\subsection{Ensembles of quadratic forms (Cor. \ref{cor quadratic}) } 
Writing each $Q_j(x)=x^TA_jx$ with $A_j$ $(d\times d)$ symmetric, the nonsingularity of $\vec{Q}$ implies that all notrivial linear combinations 
of the $A_j$ are nonsingular, which implies that at all $x-y\ne 0$ and  $\vt\in\R^k$ with $t_j\ne0$ for all $j$,
the gradients of the $Q_j(x-y)-t_j$ are linearly independent, so that $\Phi(x,y)=\vec{Q}(x-y)$ does in fact satisfy \eqref{cond dPhi}.
Furthermore, $\Sigma_{\vt}:=\{x\in\R^d: Q_1(x)=\dots=Q_k(x)=0\}$ is a
 smooth codimension $k$ surface, and the operator $\rvt$ 
is convolution with a smooth multiple of surface measure on $\Sigma_{\vt}$.

\subsection{Heisenberg spheres (Cors. \ref{cor heis} and \ref{cor step two}) } 
For simplicity, we only treat the case of the Heisenberg group $\H^1$, with  the proof for general nondegenerate step two groups being similar.
Rather than the original $\Phi:\H^1\times\H^1\to \R^2$ in \eqref{eqn orig phi for heis}, 
for simplicity we square the first component and work with
\be
\label{eqn new phi for heis}
\Phi(x,y)=\vt= \left(|x'-y'|^2, x_3-y_3+\frac12\left(x_1y_2-x_2y_1\right)\right) \in\R^2. 
\ee
Evaluating $D\Phi'=[D_x\Phi,-D_y\Phi]$, one sees that $D_x\Phi$ and $D_y\Phi$ have rank 2 
away from the lower dimensional set $Z^{sing}:=\{x'=y'\}$,
so that for $\vt\in\R_+\times\R$, $\zvt:=\{(x,y):\Phi(x,y)=\vt\}$ is smooth, codimension $k=2$ in $\H^1\times\H^1$.
On $\zvt$ we can use $x\in\H^1$ and $ \omega=(x'-y')/\sqrt{t_1}\in\sone$ as coordinates.
The phase function for $\rvt$ is
$$\phi_{\vt}(x,y,\theta)=(|x'-y'|^2-t_1)\theta_1+\left(x_3-y_3+\frac12\left(x_1y_2-x_2y_1\right)\right)\theta_2,$$
and on the resulting 6-dimensional canonical relation $C_{\vt}$ we can use $x,\, \omega$ and \linebreak$\theta\in\R^2\mzs$ as coordinates.
With respect to these,
$$\pi_L(x,\omega,\theta)=\big(x; 2\sqrt{t_1}\theta_1\omega+(1/2)\theta_2J(\omega),\theta_2\big),$$
where the cotangent variables have been split into $\xi=(\xi',\xi_3)$ and $J$ is the standard $2\times 2$ symplectic matrix 
above \eqref{eqn orig phi for heis}. From this we see that 
$$rk(D\pi_L)=4+rk\Big(\frac{D\xi'}{D(\omega,\theta_1)}\Big),$$
and this $=6$ where $\theta_1\ne 0$.
Furthermore, at the hypersurface $\{\theta_1\}\subset C_{\vt}$, 
the kernel of $D\pi_L$ is spanned by a vector with a nonzero $\partial/\partial\theta_1$ component,
so that $ker(D\pi_L)$ is transverse to $\{\theta_1=0\}$ and thus $\pi_L$ has a Whitney fold singularity
at the points where it is not a local canonical graph. By symmetry, the same is true for $\pi_R$, 
and thus $C_{\vt}$ is a folding canonical relation, so that by Thm. \ref{thm mt}, there is a loss of $1/6$ derivatives on $L^2$-based Sobolev spaces.
Hence, Thm. \ref{thm asymm} applies with $d_1=3,\, k=2$ and $\bp=1/6$. Hence, for $E,\, F\subset \H^1$ with $\hd(E)+\hd(F)>3+2+2\cdot(1/6)=16/3$, $\Delta_\Phi(E,F)$ has nonempty interior.

\section{Final comment}\label{sec comment}

We observe that the threshold in Thm. \ref{thm same d} cannot  in general be lowered.
Returning to the setting of Steinhaus' theorem, one has $k=d$,  $\Phi(x,y)=x-y$, and for each $\vt\in\R^d$ the operator $\rvt$ is just translation by $\vt$,
which is an FIO  of order 0 associated with a canonical graph.
Hence $\bp=0$ and $\frac{d+k}2+\bp=d$. 
The resulting sufficient lower bound in Thm. \ref{thm same d} is then  $\hd(E)>d$,
so that the theorem is vacuous in this case, and   does not imply Steinhaus' result.
However,  the threshold $\frac{d+k}2+\bp$ cannot be lowered in this case: 
Falconer's original counterexamples related to the distance problem (see \cite[Thm. 2.4]{Falc86})  
can be modified to show that for any $s<d$ there is an $E\subset \R^d$ with $\hd(E)=s$ and $\hbox{int}(E-E)=\emptyset$,
showing that the range of $\hd(E)$ in Thm. \ref{thm same d} cannot {\it in general} be lowered below the endpoint 
$\frac{d+k}2+\bp$.
Of course, this leaves open the possibility of  improvement for other, specific $\Phi$.


\end{document}